\documentclass[a4paper,11pt]{article}
\usepackage{graphicx,amsmath, amssymb,a4wide,amsthm}
\usepackage[numbers,sort&compress]{natbib}
\usepackage{bm}
\usepackage[english]{babel}
\usepackage{hyperref}
\hypersetup{
    bookmarksopen=false,
    bookmarksnumbered=true,
    pdftitle={Applications of polling systems},
    pdfauthor={M.A.A. Boon, R.D. van der Mei, E.M.M. Winands}
}
\ifpdf
  \hypersetup{colorlinks=true,linkcolor=black,urlcolor=black,citecolor=black,pdfpagemode=UseOutlines,plainpages=false,pdfpagelabels}
\else
  \hypersetup{colorlinks=false}
\fi

\newtheorem{theorem}{Theorem}[section]

\newtheorem{property}[theorem]{Property}

\theoremstyle{definition}

\providecommand{\href}[2]{#2}

\newcommand{\z}{z_1,\dots,z_N}

\newcommand{\E}{\mathbb{E}}

\renewcommand{\L}{\widetilde L}
\newcommand{\B}{\widetilde{B}}
\newcommand{\equaldist}{\,{\buildrel d \over =}\,}
\renewcommand{\S}{\widetilde{S}}
\newcommand{\W}{\widetilde{W}}

\title{Applications of polling systems}

\author{M.A.A. Boon\footnote{\textsc{Eurandom} and Department of Mathematics and Computer Science, Eindhoven University of Technology, P.O. Box 513, 5600MB Eindhoven, The Netherlands}\\\href{mailto:marko@win.tue.nl}{marko@win.tue.nl} \and R.D. van der Mei\ \footnote{Department of Mathematics, Section Stochastics, VU University, De Boelelaan 1081a, 1081HV Amsterdam, The Netherlands}\ \ \footnote{Centre for Mathematics and Computer Science (CWI), Department of Probability and Stochastic Networks, 1098 SJ Amsterdam, The Netherlands}\\\href{mailto:mei@cwi.nl}{mei@cwi.nl} \and E.M.M. Winands \footnotemark[2]\\\href{mailto:emm.winands@few.vu.nl}{emm.winands@few.vu.nl}}

\date{February, 2011}

\begin{document}
\maketitle

\begin{abstract}
Since the first paper on polling systems, written by Mack in 1957, a huge number of papers on this topic has been written. A typical polling system consists of a number of queues, attended by a single server. In several surveys, the most notable ones written by Takagi, detailed and comprehensive descriptions of the mathematical analysis of polling systems are provided. The goal of the present survey paper is to complement these papers by putting the emphasis on \emph{applications} of polling models. We discuss not only the capabilities, but also the limitations of polling models in representing various applications.
The present survey is directed at both academicians and practitioners.

\end{abstract}

\section{Introduction}

A typical \textit{polling system} consists of a number of queues, attended by a single server in a fixed order. There is a huge body of literature on polling systems that has developed since the late 1950s, when the papers of Mack et al. \cite{mack2,mack1} concerning a patrolling repairman model for the British cotton industry were published. The term ``polling'' originates from the so-called polling data link control scheme, in which a central computer (server) interrogates each terminal (queue) on a communication line to find whether it has any information (customers) to transmit. The addressed terminal transmits information and the computer then switches to the next terminal to check whether that terminal has any information to transmit. In a broader perspective, polling models are applicable in situations in which several types of users compete for access to a common resource which is available to only one type of user at a time. The ubiquity of polling systems can be observed in many applications, e.g., in computer-communication, production, transportation and maintenance systems. The present paper discusses the main application areas of polling systems along with an extensive list of references, examines how these various applications can be represented and analyzed via polling models and outlines several topics for future research inside and outside the described areas.

To illustrate the fact that polling systems have become a powerful tool for the performance analysis of a wide variety of important applications, we would like to cite Hideaki Takagi, one of the fathers of the success of polling models \cite{takagi3}: \textit{The analysis of polling models gained momentum as queueing systems that are easy to understand, analyze, and extend. The study has been accelerated largely by applications to the modeling of communication, manufacturing and transportation systems. I believe that it is one of the few successful theoretical performance evaluation models developed in the last decades.} The fact that nowadays polling systems are still fully alive can be illustrated by observing that the recent application-oriented conferences Performance 2007, Informs Applied Probability 2007 and 2009, and ValueTools 2008 all scheduled a dedicated session for polling systems. Very recently, Annals of Operations Research dedicated a special issue to polling systems.

Until the mid nineties, Takagi maintained a fairly complete
bibliography on polling models, which contained over 700 publications, including journal and conference papers, books,
theses, and technical reports. Since research on polling models has continued in the past years
(although not at as high a pace as in the years around 1990), the number of papers might now well be above 1000. The main reason for this tidal wave of papers on polling models  in the past 50 years can be found in the  diversity of applications, which have led to innumerable enhancements of the basic cyclical polling system. Such enhancements gave rise to interesting new queueing theoretic problems which in turn stimulated research in various directions.

In this respect one of the most remarkable mathematical results
in the polling literature is the striking dichotomy in the complexity of the analysis of polling
systems, which has been independently illuminated by Fuhrmann \cite{fuhrmann} and Resing \cite{resing} via the use of a multitype
branching approach. That is, if the service discipline (i.e., the criterion which determines how many customers are served during a visit of the server to a queue)
satisfies a certain branching property, in many cases the polling system allows for an exact analysis by rather standard
methods. If this branching property is, however, violated, the corresponding polling system
defies an exact analysis except for some special (two-queue or symmetric) cases.

The most important members of the class of policies satisfying the branching property are the
\textit{exhaustive} service discipline, under which the server continues to work until the
queue has become empty, the \textit{gated} service discipline, under which exactly
those customers are served who were present at the queue at the beginning of
the visit, and \textit{globally gated} service, which states that only those customers are served that were present at the server's arrival epoch at some predetermined ``parent queue''. (Actually the globally gated service discipline does not strictly satisfy the branching property as defined in \cite{resing}, but a weaker variant that still allows an exact analysis.) Unfortunately, the \textit{$k$-limited} service policy, under which the server continues working at a queue until either a predefined number of $k$ customers is served or until the queue becomes empty, is on the wrong side of the borderline implying that even mean queue lengths are in general not known.

In the past various survey papers dedicated to polling systems have appeared in the open literature (see, e.g., \cite{levy1,takagi1,takagi2,takagi3,yechialisurvey93,vishnevskii}) which give detailed and comprehensive descriptions of the
mathematical analysis of polling systems. Although we present a brief overview of the basic mathematical analysis techniques, the main goal of the present survey paper is to complement these papers by putting the emphasis on \emph{applications} of polling models. We unearth not only the capabilities, but also the limitations of polling models in representing various applications. Moreover, since the publication of the previous survey papers around 10 till 15 years ago, a large number of papers on polling models have been published. The present survey is directed at both academicians and practitioners. The list of references is relatively recent, but deliberately far from exhaustive. The reader who is interested in going deeper into specific mathematical details of polling systems will find references to the aforementioned survey papers aimed at that objective. In selecting the literature to cite in this survey, a preference has been given to current papers and research on applications of polling systems.

The rest of the present paper is organized as follows. Section \ref{pollingsection} gives a global overview of various aspects of polling systems, and gives an illustration of how one can analyze a typical polling system. In Section 3 we describe the main applications of polling systems in communication networks, production systems, traffic and transportation problems and various other fields of engineering. Finally, the last section indicates some possible directions for further research.

\section{Polling systems}\label{pollingsection}

A polling model consists of a number of queues, attended by a single server who visits the queues in some order to render service to the customers waiting at the queues, typically incurring some switch-over time while moving from one queue to the next. In Section \ref{features} we give a structural overview of the variety of polling models considered in the literature, successively discussing model variants with respect to the arrival process, the buffer size, the service process,
the switch-over process, the server routing, the service discipline and the queueing
discipline. In Section \ref{mathematicalanalysis} we briefly illustrate how polling systems may be analyzed.

\subsection{Typical features}\label{features}

\paragraph{Arrival process.} In most polling models, customers arrive at
the queues according to mutually independent homogeneous single Poisson arrival
processes. In many cases arrival processes can indeed be described adequately
by Poisson processes, e.g., in the case of telephone calls and traffic accidents.
Moreover, due to the memorylessness property of the Poisson process, polling
models with Poisson arrivals often allow for a tractable analysis.

Yet, in many situations the assumption of Poisson arrivals is unrealistic, e.g., in case of bursty traffic (voice, video) and several production applications. Only very recently, variations of Poisson arrival processes have been considered. Van Vuuren and Winands \cite{vuuren} have presented an approximate algorithm for polling systems with $k$-limited service and Markovian Arrival Processes (MAPs). Boxma et al. \cite{pollinglevy09} consider polling systems with L\'evy-driven input. Bertsimas and Mourtzinou \cite{bertsimas99} study polling systems with Mixed General Erlang arrival processes. Slightly more general arrival processes are discussed by Saffer and Telek \cite{saffertelek2010} who consider BMAPs (Batch Markov Arrival Processes). Some limiting cases lead to closed-form expressions for, e.g., queue lengths or waiting time distributions, even under the assumption of general renewal arrivals. For example, Van der Mei and Winands \cite{mei1,vdmeiwinands08} derive the distribution of the scaled waiting time under heavy traffic conditions, i.e., when the load tends to one. For switch-over times tending to infinity, the scaled waiting time distribution is conjectured in \cite{winandslargesetupsbranching09}, but the proof is still an open problem. Closed-form approximations for waiting times in polling systems with general renewal arrivals are constructed in \cite{boonapprox2011,dorsman2010} using an interpolation between light-traffic and heavy-traffic limits.

In polling models it is almost exclusively assumed that customers arrive from some external infinite population. In some applications this assumption is unrealistic, e.g., in maintenance environments where customers represent machine breakdowns. Altman and Yechiali \cite{altman2} study a closed polling model in which customers, after having received service at a queue, proceed to another (possibly the same) queue. This model is extended in \cite{armonyyechiali99} where a model with transient as well as permanent customers is considered. Sidi and Levy \cite{sidi} consider a system in which (external) customers arrive at a queue according to a Poisson process, and in which a customer, after having received service at a queue, either leaves the system or moves to another queue (with some given probability). Applications of this type of customer routing can be found, for example, in manufacturing when products undergo service in a number of stages or in the context of rework. Finally, we mention work of Levy and Sidi \cite{levy5} in which they study a polling model with simultaneous arrivals at the queues.

A final remark regarding arrival processes, is that some papers study discrete-time arrival processes. For some applications, it is natural to divide time into slots. In communication systems discrete-time models are used by, e.g., Kleinrock and Scholl \cite{kleinrockscholl} for the analysis of the Minislotted Alternating Priorities (MSAP) multiple-access scheme, by Kleinrock and Levy \cite{kleinrock} for the analysis of polling systems with random server routing to model the Slotted ALOHA system, and by Beekhuizen \cite{Beekhuizen2010} to model networks on chips. An example in a completely different application area, is provided by Van Leeuwaarden \cite{leeuwaarden06} who analyses the Fixed Cycle Traffic Light queue in discrete time.

\paragraph{Buffer size.} In most polling models the buffer size is assumed to be infinite. In a number of applications the buffer capacities are obviously finite, e.g., transportation and manufacturing. In some applications it is natural to assume that the buffers are unit sized, i.e., each queue can only accommodate one customer at a time. Applications of polling models with unit-sized buffers are found, e.g., in the area of maintenance. If one would like to study a multi-item production system which produces products in a make-to-order fashion, a polling system with finite buffers appears naturally as well. For an analysis of polling models with finite buffers, see, e.g., Takagi \cite{takagi4}, or Browne and Yechiali \cite{browneyechiali}.

\paragraph{Service process.} The service times at a queue are typically assumed to be samples from a probability
distribution which is characteristic for that queue. The service times
are usually assumed to be mutually independent and independent of the actual
state of the system. Polling models with batch service are studied in \cite{vdwalyechiali2003} and \cite{boxmavdwalyechiali08}. Van der Wal and Yechiali \cite{vdwalyechiali2003} study a model where, at the beginning of each cycle, the server determines a tour along each of the non-empty queues. This model, with customers being served in one batch, has applications to videotex, telex, and  time-division multiple access (TDMA) systems, as well as for central database systems and Automated Guided Vehicles (AGVs). Boxma et al. \cite{boxmavdwalyechiali08} study a polling model with batch service. Server routing is either strictly cyclic, or the server visits the queue with the oldest customer first.

The usual assumption in the polling literature is that a \textit{single} server is visiting all the queues. Only few results have been obtained on \textit{multiple}-server polling systems. For example, some numerical results are obtained by Ajmone Marsan et al. \cite{ajmonemarsan94}, who perform several numerical studies using generalized stochastic Petri nets (GSPNs). Analytical approximative results have been obtained by Morris and Wang \cite{morriswang84}, and also in \cite{borst3,mei0}. Morris and Wang observe the interesting phenomenon that the servers tend to cluster if they follow identical routes, deteriorating the system performance. Because the methods used in all of these papers are based on approximations, many challenges remain for future research. However, some interesting exact results are found for a polling system with multiple \emph{coupled} servers, visiting the queues always together \cite{borst2}.

\paragraph{Switch-over process.}
The switch-over times needed by the server to move from one queue to another
queue are typically assumed to be samples from some prespecified probability
distribution which is characteristic for that couple of queues. In a few cases the
switch-over times are assumed to be decomposed into switch-out times (which
are characteristic for the queue being departed from) and switch-in times
(which are characteristic for the queue being switched to), putting a special
structure on the switch-over times.

Nearly all of the existing literature on polling systems makes the assumption of state-independent switch-overs, i.e., switch-overs are assumed to be independent of the current state of the system. Notable exceptions are the recent studies of Altman et al. \cite{altman}, G\"unalay and Gupta \cite{gunalay}, Gupta and  Srinivasan \cite{gupta}, Singh and Srinivasan \cite{singh} and Winands et al. \cite{winands2}. The choice of modeling state-independent setups is generally not motivated by an application but by the tractability of the resulting analysis.

Altman and Fiems \cite{altman4} allow correlation between switch-over times in the polling model. That is, they assume that the switch-over times constitute a stationary ergodic series of random variables. A wireless LAN, where an access point polls mobiles, is one of the applications of this type of switch-over process.

\paragraph{Server routing.}
The order in which the server visits the queues is determined by some routing
mechanism. Such a mechanism may depend on the actual state of the system
(dynamic) or may be independent of the state of the system (static).

We first discuss \textit{static} routing mechanisms. The traditional routing mechanism
is the cyclic server routing. To model systems in which particular queues are
visited more frequently than others, cyclic polling has been extended to periodic
polling, in which the server visits the queue periodically according to some
service order table of finite length, cf. \cite{baker}. Alternatively, the server may
be routed along the queues according to some probabilistic routing mechanism.
Kleinrock and Levy \cite{kleinrock} introduce the so-called random polling mechanism
in which, after a departure from a queue, the server proceeds to queue $j$ with
some given probability $p_j$.

Boxma and Weststrate \cite{boxma} and Srinivasan \cite{srinivasan91} generalize the random
polling scheme to the so-called Markovian routing scheme in which the server,
after a departure of the server from queue $i$, proceeds to queue $j$ with given
probability $p_{i,j}$.

Under \textit{dynamic} server routing, the decision of the server as to the order in
which the queues are visited may depend on a certain amount of information
available to the server, such as the queue lengths. For instance, it might not
make sense to move to an empty queue while customers are waiting at other
queues (cf., e.g., \cite{hofriross87}). An example of dynamic server routing in systems with full information
about the buffer contents is the 'serve-longest-queue' policy.

Yechiali \cite{yechiali} has introduced the so-called \textit{semi-dynamic} server routing,
in which the server, at the end of each tour along the queues, receives
information about the queue lengths. Based on this information, the server
makes a decision as to the order in which it will visit the queues in the
next tour. This order cannot be changed during the course of that visit tour.

\paragraph{Service discipline.}
The service discipline specifies the number of customers that is served during
one visit of the server to a queue. 
In the literature, a whole abundance of service disciplines has been proposed by
putting some cut-off mechanism (which eventually depends on the evolution of
the queue length during the server visit) on the classical exhaustive and gated
service disciplines. The service disciplines can be classified into the class of
\textit{customer-limited} service disciplines, in which restrictions are put on the number
of customers served during a visit of the server to a queue, and the class
of \textit{time-limited} service strategies, putting restrictions on the amount of time
spent by the server during one visit of the server to a queue.

Alternatively, service disciplines can be classified into the so-called \textit{exhaustive-type} policies and the \textit{gated-type} policies, depending on whether customers who arrive at a queue while the server is working at that queue are candidates for service during the same visit of the server to that queue, e.g., \cite{levy1}. Under an exhaustive-type policy the customers arriving at a queue in service are \textit{candidates} for service in the same visit period, whereas under a gated-type policy they are not. The best known variant of the gated service discipline is \emph{globally gated}, which also allows for an exact analysis \cite{boxmalevyyechiali92}. This service discipline states that only those customers are served who were present at the beginning of a visit to one fixed ``parent-queue''. Numerous hybrid variants of service disciplines can be conceived by combining the various types of cut-off mechanisms. A number of these service disciplines has been studied in the literature, such as customer-limited type service policies like the $k$-limited service, Bernoulli service~\cite{keilson2}, probabilistically-limited service \cite{leung}, binomial-gated service \cite{levy4}, binomial-exhaustive, Bernoulli-gated and Bernoulli-exhaustive service \cite{browneyechiali2}, fractional-exhaustive service \cite{levy3}, mixed gated/exhaustive service based on customer priority levels \cite{boonadan2008}, Bernoulli-type service \cite{resing}, and time-limited service disciplines with exponential time limits \cite{coffman,leung2}, and with constant time limits \cite{souza}. A gated-type service discipline with multiple gates is studied in \cite{meiresing07,meiresing08,meiroubos2011}.

\paragraph{Queueing discipline.}
The queueing discipline specifies the order in which the customers present at
the same queue are served. The most common queueing discipline is the classical
\textit{First-Come-First-Served} (FCFS) discipline. Very recently, Wierman et al.
\cite{wierman} showed that optimization of the queueing discipline in a polling system can improve the performance of the system significantly without
having to purchase additional resources. Other recent papers in this direction
are, for example, \cite{boon,boonpriorities,boxma6}. It is important to note that
the queue length distribution is independent of the service order, and so are,
by Little's law, the mean waiting times (provided the queueing discipline does
not depend on the service times). However, the distribution of the waiting time
does depend on the queueing discipline.

\subsection{Mathematical analysis}\label{mathematicalanalysis}

In this subsection we give a brief impression of how polling systems can be analyzed. As stated in the introduction, the present paper surveys the application areas of polling systems, where it is not our intention to present an encyclopedic overview of all available mathematical results. Consequently, we have chosen to discuss only one (frequently used) analysis method, studying the marginal queue length and waiting time distributions. The reader should bear in mind that several alternatives exist. For detailed and comprehensive surveys of the analysis methodology of polling systems, we refer to Takagi's surveys \cite{takagi1,takagi2,takagi3}, and furthermore to \cite{levy1,yechialisurvey93,vishnevskii}. Moreover, Takagi's seminal book ``Analysis of polling systems'' \cite{takagi86} should definitely be mentioned here.

The model studied in the present subsection can be considered as the basic, standard polling model. We assume that the polling model consists of $N$ queues $(N\geq 1)$, denoted by $Q_1,\dots,Q_N$, attended by one server in cyclic order. Whenever the server switches from $Q_i$ to $Q_{i+1}$, a switch-over time $S_i$ is incurred ($i=1,\dots, N$; since the system is cyclic, $Q_{N+1}$ actually refers to $Q_1$ etc.). The switch-over times are independent random variables. 
In \cite{Borst1} an adaptation of the analysis is given for polling models without switch-over times. The arrival processes to the different queues are independent Poisson processes with intensities $\lambda_i$, $i=1,\dots,N$. The service times of customers arriving in $Q_i$ are independent random variables denoted by $B_i$. At the end of this section we will briefly mention some alternative techniques that may be required if the polling model of interest does not satisfy all of these assumptions.

The individual loads of the queues are denoted by $\rho_i=\lambda_i\E[B_i]$, for $i=1,\dots,N$, and the total load of the system is denoted by $\rho=\sum_{i=1}^N \rho_i$. The Laplace-Stieltjes transform (LST) of a continuous, positive random variable $X$ is denoted by $\widetilde{X}(\cdot)$. We use the same notation for the probability generating function (PGF) in case $X$ is a discrete random variable.

\paragraph{Cycle times, visit times and intervisit times.}
The distribution of the length of one cycle depends on its starting point. For this reason we denote by $C_i$ the time between two successive visit \emph{beginnings} to $Q_i$, and by $C_i^*$ the time between two consecutive visit \emph{completions} to $Q_i$. The \emph{intervisit time} of $Q_i$, denoted by $I_i$, is the time between the visit completion of $Q_i$ and the next visit beginning to this queue. The \emph{mean} cycle time and the \emph{mean} visit times can easily be computed in a cyclic polling system. The first obvious, but important, observation is that the mean cycle time does \emph{not} depend on the starting point of the cycle. For this reason, we denote the mean cycle time by $\E[C]=\E[C_i]=\E[C^*_i]$ for all $i=1,\dots,N$. Since each customer arriving in the system gets served eventually, the server must be working a fraction $\rho$ of the time, and is switching between queues a fraction $1-\rho$ of the time.  This leads to the simple relation $\E[S]=(1-\rho)\E[C]$, where $S=\sum_{i=1}^N S_i$. A similar way of reasoning can be applied to the workload of the individual queues, yielding expressions for the mean visit and intervisit times.
\begin{align*}
\E[C] &= \frac{\E[S]}{1-\rho},\\
\E[V_i] &= \rho_i\E[C],\\
\E[I_i] &= (1-\rho_i)\E[C].
\end{align*}
The stability condition of the system depends on the service disciplines in the various queues. For exhaustive, gated, and globally gated service (and many other disciplines as well) the stability condition is simply $\rho<1$, yielding $\E[C]<\infty$. Note that for these service disciplines the switch-over times do not play a role in the stability condition, as they become negligible compared to the visit times, as $\rho \rightarrow 1$.

\paragraph{Queue lengths and waiting times.}
The key performance measures of a polling system are the queue lengths and the waiting times. One of the most remarkable results regarding the marginal queue length distribution, denoted by $L_i$, has been obtained by Fuhrmann and Cooper \cite{fuhrmanncooper85}. They observed that the marginal queue length distribution can be decomposed as the sum of two independent random variables:
\[
L_i \equaldist L_i^{(M/G/1)} + L_i^{(I_i)},
\]
where $L_i^{(M/G/1)}$ is the queue length of an $M/G/1$ queue with only type $i$ customers in isolation, often referred to as ``the corresponding $M/G/1$ queue'', and $L_i^{(I_i)}$ is the queue length at an arbitrary moment during the intervisit period $I_i$. This results in the following PGF of the marginal queue length distribution:
\begin{equation}
\L_i(z) =\,\frac{(1-\rho_i)(1-z)\B_i\left(\lambda_i(1-z)\right)}{\B_i\left(\lambda_i(1-z)\right)-z}
\times\frac{\L_i^{(V_{c_i})}(z) - \L_i^{(V_{b_i})}(z)}{\lambda_i\E[I_i](1-z)},\label{fuhrmanncooperQLdecomposition}
\end{equation}
for $i=1,\dots,N$. In \eqref{fuhrmanncooperQLdecomposition} we have used $L_i^{(V_{b_i})}$ and $L_i^{(V_{c_i})}$ to denote the number of customers in $Q_i$ at, respectively, the beginning and completion of a visit to $Q_i$. An important observation is that one only needs the queue length distribution at the beginnings and completions of a visit to each queue, in order to find the marginal queue length distributions. Moreover, Keilson and Servi show that the LSTs of the waiting time distributions follow directly from $\L_i(\cdot)$ by applying the distributional form of Little's Law \cite{keilsonservi90}:
\[
\W_i(\omega) = \frac{\L_i(1-\omega/\lambda_i)}{\B_i(\omega)},
\]
where $\W_i(\cdot)$ is the LST of the waiting time distribution of type $i$ customers, \emph{excluding} the service time $B_i$.
\paragraph{Joint queue length at polling epochs.}
Determining the marginal queue length distributions at visit beginnings and completions, requires studying the \emph{joint} queue length distributions at these embedded epochs. It turns out that, a few exceptions aside, an exact analysis is only possible if all service disciplines in the polling system satisfy the following \emph{branching property}, described by Fuhrmann \cite{fuhrmann} and Resing \cite{resing}.
\begin{property}[Branching Property]\label{resingproperty}
If the server arrives at $Q_i$ to find $k_i$ customers there, then during the course of the server's visit, each of these $k_i$ customers will effectively be replaced in an i.i.d. manner by a random population having PGF $h_i(z_1,\dots,z_N)$, which can be any $N$-dimensional PGF.
\end{property}
The PGF mentioned in the branching property depends on the service discipline. As an example, we give the expressions for exhaustive and gated service.
\begin{align*}
&\textrm{Gated: } && h_i(z_1,\dots,z_N) = \B_i\big(\sum_{j=1}^N\lambda_j(1-z_j)\big),\\
&\textrm{Exhaustive: } && h_i(z_1,\dots,z_N) = \pi_i\big(\sum_{j\neq i} \lambda_j(1-z_j)\big),
\end{align*}
where $\pi_i(\cdot)$ is the LST of a busy period distribution in an $M/G/1$ system with only type $i$ customers, so it is the root in $(0,1]$ of the equation $\pi_i(\omega) = \B_i\big(\omega + \lambda_i(1 - \pi_i(\omega))\big)$, $\omega \geq 0$ (cf. \cite{cohen82}, p. 250).

Recall that the globally gated service discipline does not strictly satisfy Property \ref{resingproperty}, but still allows for an exact analysis. Borst \cite{borst} defines a generalization of the branching property that still ensures mathematical tractability. The globally gated discipline satisfies this weaker property.

Assuming that each queue satisfies the branching property, the PGF of the joint queue length distribution at visit completions, denoted by $\L^{(V_{c_i})}(\z)$, can be expressed in terms of the joint queue length distribution at visit beginnings, denoted by $\L^{(V_{b_i})}(\z)$. The following relations, which are commonly referred to as the \emph{laws of motion}, hold:
\begin{align}
\L^{(V_{c_1})}(\z) &=\, \L^{(V_{b_1})}\big(h_1(\z), z_2, \dots, z_N\big),\label{lawsofmotion1}\\
\L^{(V_{b_2})}(\z) &=\, \L^{(V_{c_1})}(\z)\ \S_1\big(\sum_{j=1}^N\lambda_j(1-z_j)\big),\\
&\ \ \vdots\nonumber\\
\L^{(V_{c_N})}(\z) &=\, \L^{(V_{b_N})}\big(z_1, z_2, \dots, z_{N-1},h_N(\z)\big),\\
\L^{(V_{b_1})}(\z) &=\, \L^{(V_{c_N})}(\z)\ \S_N\big(\sum_{j=1}^N\lambda_j(1-z_j)\big).\label{lawsofmotionN}
\end{align}
The recursive equation for $\L^{(V_{b_1})}(\z)$ itself can be used to compute the moments of this joint queue length distribution. An explicit expression for $\L^{(V_{b_1})}(\z)$ can be found by clever definition of \emph{offspring generating functions} and \emph{immigration generating functions}. We refer to \cite{resing} for more details.

\paragraph{Mean waiting times.}
By differentiation of the LSTs of the waiting time distributions, the following expressions for $\E[W_i]$ are found for exhaustive, gated, and globally gated service respectively:
\begin{align}
&\textrm{Exhaustive: } & \E[W_i] &=\, (1-\rho_i)\frac{\E[{C^*_{i}}^2]}{2\E[C]},   \label{EWiExhaustive}\\
&\textrm{Gated: } & \E[W_i] &=\, (1+\rho_i)\frac{\E[C_{i}^2]}{2\E[C]},\label{EWiGated}\\
&\textrm{Globally gated: } & \E[W_i] &=\, \sum_{j=1}^{i-1}\E[S_j]+\big(1+2\sum_{j=1}^{i-1}\rho_j + \rho_i\big)\frac{\E[C_{1}^2]}{2\E[C]}.\label{EWiGloballyGated}
\end{align}
For a polling system with globally gated service, we have assumed that the server's arrival epoch at $Q_1$ determines which customers (in \emph{all} queues) are eligible for service. Although the expressions \eqref{EWiExhaustive}--\eqref{EWiGloballyGated} seem quite simple at first sight, it should be noted that the computation of the second moments of $C_i^*$ and $C_i$ is not straightforward. Although they can be obtained by solving the set of equations resulting from differentiation of the recursive expression for $\L^{(V_{b_1})}(\z)$, obtained from the laws of motion, with respect to $z_1, \dots, z_N$, more efficient methods have been developed. The two most frequently used alternatives to find mean waiting times are the \emph{descendant set approach}, developed by Konheim et al. \cite{konheim94}, and mean value analysis for polling systems, developed by Winands et al. \cite{winands06}.

\paragraph{Pseudo-conservation law.} Boxma and Groenendijk~\cite{boxmagroenendijk87} derive a so-called \emph{pseudo-conservation law} (PCL) for the mean waiting times under various service disciplines. This PCL states that
\begin{align}
\sum_{i=1}^N \rho_i \E[W_i]&=\rho \frac{\sum_{i=1}^N \lambda_i \E[B_{i}^2]}{2(1-\rho)}+\rho \frac{\E[S^2]}{2\E[S]}
 +\frac{\E[S]}{2(1-\rho)} \left( \rho^2 - \sum_{i=1}^N \rho_i^2\right)+\sum_{i=1}^N\E[Z_{ii}],\label{pcl}
\end{align}
where $Z_{ii}$ denotes the amount of work left behind by the server at $Q_i$ at the completion of a visit, which solely depends on the service discipline of $Q_i$. For exhaustive service at $Q_i$, we have $\E[Z_{ii}]=0$, and for gated service $\E[Z_{ii}]={\rho_i^2\E[S]}/({1-\rho})$.

\paragraph{Other models and/or techniques}
In the present section we have described only one technique, often referred to as the \emph{buffer occupancy method}, first used by Cooper and Murray \cite{coopermurray69} and Cooper \cite{cooper70}. An advantage of this technique is its applicability to several variations of the basic polling model, described in the current subsection, such as polling systems with other service disciplines than exhaustive and/or (globally) gated service, systems with simultaneous arrivals \cite{levy5}, or fluid polling systems \cite{Czerniak2009}. However, if the assumptions stated in the beginning of this subsection are not all valid, one might have to resort to alternative ways to analyze the system. For example, if the arrival processes are not Poisson, but BMAPs, a generalization of the buffer occupancy method using Kronecker product notation can be used \cite{saffertelek2010}. Altman and Fiems \cite{altman4} show how stochastic recursive equations can be used to analyse a polling model with correlated switch-over times. If the service disciplines do not satisfy the branching property, one may have to resort to different techniques, or even approximations. For example, a two-queue polling system with Bernoulli service (a generalization of 1-limited as well as exhaustive service) has been analyzed by Lee \cite{Lee} by formulating the model as a Riemann boundary value problem with a shift. Unfortunately, to analyze polling systems with non-branching service disciplines consisting of more than two queues, one may have to use approximations or numerical methods. Van Houdt \cite{vanhoudt2010} uses a numerical algorithm based on the power method, Kronecker matrix representations, and the shuffle algorithm, to analyze discrete-time polling systems for networks on chips (which will be discussed later in this paper). This algorithm can also be used for general polling systems with branching-type or non-branching-type service disciplines.

\section{Applications}

In this section, we will first describe the three most successful application areas of polling models
within the classical fields of engineering of communication systems, production systems, traffic and transportation systems. Subsequently, the last subsection summarizes a variety of miscellaneous applications of polling systems.

\subsection{Computer-communication systems}

Polling systems find a wealth of applications in the area of computer-communication systems, where
resources (e.g., bandwidth, CPU power) are shared among different users. The reader is referred to surveys by Grillo
\cite{Grillo1}, Levy and Sidi \cite{levy1}, Takagi \cite{takagi3} and Weststrate \cite{weststrate} for overviews
of applications up until the early 1990s. For completeness, the main applications cited in these survey papers are
outlined below, and supplemented with more recent applications of polling models in communication systems.

\paragraph{Time-sharing computer systems.}
Classical applications of polling models are time-sharing computer
systems \cite{Klimov1}, consisting of a
number of terminals connected by multi-drop lines to a central computer. The data transfer from the terminals
to the computer - and back - is controlled via a polling scheme in which the computer polls the terminals,
requesting their data, one terminal at a time. In such applications of polling models, the server represents
the central computer, the queues represent the terminals and customers represent the data.

In communication networks different terminals compete for access to a shared medium. If multiple terminals transmit
or receive data over the medium at the same time, packet collisions and interference problems may occur. Motivated by
this many medium access control (MAC) protocols have been proposed for different network technologies, in many cases
leading to polling models.

\paragraph{Token-ring networks.}
Bux \cite{Bux1} uses polling models to study the performance of token-passing schemes in Local Area Networks (LANs) where a token - representing the right for transmission - is circulated among the different users. In such cases, the token-passing scheme is usually configured in a ring or a bus topology. A token-ring network can be characterized as a set of stations connected to a common transmission medium in a ring topology. All messages travel over a fixed route from station to station around the loop. A token can be in two states: occupied or free. A station with data to transmit reads the free token and changes it to the occupied state before retransmitting it. The occupied token is then incorporated as part of the header of the data transmitted on the ring by the station. Thus, other stations on the ring can read the header, note the occupied token and refrain from transmission. When the token is back at the station that changed it to the occupied state and the station decides to transfer the right for transmission to another station, it changes the token to the free state. The token-ring network allows the transmission of packets in a conflict-free manner. However, transmission of packets may still fail due to errors and distortions on the ring itself. Typically, these errors are rare and a so-called Selective-Repeat ARQ (SR-ARQ) could be used to recover these errors. In such a scheme, a station that receives an erroneous message transmits a negative acknowledgment to the transmitting station to indicate that the message has to be retransmitted. To analyze the performance of SR-ARQ schemes, Levy and Sidi \cite{levy1} use a polling model in which each station is represented by two queues, one for messages that need to be sent out and one for negative acknowledgments to be sent back when erroneous messages are received. Altman and Kofman \cite{altmankofman} propose a solution to deal with the irregular, bursty, correlated arrival processes in token-ring networks. To this end, they use a polling model with so-called Cruz-type traffic, filtering the arrival streams by leaky buckets.

\paragraph{Token-bus networks.}
The token-bus network consists of a set of stations connected to each other in a bus topology. The intention behind this technique is to combine the attractive features of the bus topology with those of a conflict-free medium-access protocol. In a token bus, as the token is passed, a logical ring is formed. Since the bus topology does not impose any sequential ordering of the stations, the logical ring is defined by a sequence of station addresses. Conceptually, token passing on buses and rings is very similar and the same type of performance model can be used to describe these techniques. The difference between a token ring and a token bus from a modeling point of view is that the server in the token ring network visits the queues in a cyclic manner, whereas in the token-bus model the server moves along the queues in a non-cyclic periodic manner, which can be modeled by a polling table. An example of a token-bus network is presented by Manfield \cite{Manfield1}, where a communication network constitutes the transmission medium between a master processor and a set of peripheral processors. The operation of these systems works as follows: the master processor polls each of the peripheral processors in turn. Subsequently, it receives a response which indicates whether the peripheral processor has any packet to send. If it has, the packet at the head of that queue is sent to the master processor. After polling a peripheral incoming queue, the master processor checks its own queue for outgoing packets, and if the queue is not empty, seizes control of the bus. The advantage of this system, is that it ensures that outgoing messages are rapidly sent and it avoids message deadlocks in the central control facility (i.e., the master processor). This polling scheme is called \emph{star polling} in the literature.

\paragraph{Slotted-ring networks.}
Another class of communication protocols for networks with a ring topology is a so-called slotted-ring network. In a slotted-ring network one or more slots circulate along the stations. If there is a packet at a station ready for transmission and an empty slot comes along, the packet is put into the slot, together with the address of the destination station. That slot is then examined by each of the other stations in turn, until the destination station recognizes it and copies its content. There are two possibilities to empty the slot: either it is emptied by the source station (this is called the Cambridge Ring Protocol), or
by the destination system (this is called the Orwell Ring Protocol). We refer to Bux \cite{Bux1} for a polling model with source
release, and to Van Arem \cite{Arem1} for a polling model of a slotted-ring protocol with destination release.


\paragraph{Fibre Distributed Data Interface networks.}
The FDDI is a token-passing protocol for LANs with a ring topology in which the
access to the ring for transmission is controlled via a so-called timed-token protocol, i.e. where the transmission time for
each station is bounded. This type of models leads to the formulation of polling models with time-limited service policies \cite{Johnson1}.

\paragraph{Distributed Queue Dual Bus networks.}
The DQDB protocol is a multiple access protocol for communication networks consisting
of two unidirectional buses carrying information in the opposite directions. The stations are distributed along the two buses and have the capability to transmit/receive information to/from both buses. The DQDB protocol is intended to integrate data, voice and video
traffic in a single communication network. Bisdikian \cite{Bisdikian1} uses a polling model to study medium access mechanisms of a single station
in a DQDB network.

\paragraph{Random access schemes.}
Opposite to the scheduled multiple access protocols are the random multiple access protocols, where an entity with a message will transmit it regardless of potential collisions. The ALOHA access scheme is an example of a random access scheme. As soon as a packet arrives at a station it is transmitted. When the transmission fails, for example because other packets were transmitted at the same time causing a collision, the packet is scheduled for transmission after a random period of time.
An alternative  for this scheme is the so-called reservation ALOHA access scheme. In this scheme, a station is granted the exclusive right to transmit, without being interfered by any other station, for a certain amount of time. When a transmitting station no longer reserves the channel, some - or all - stations of the system start contending in order to seize the channel. The length of the contention period is random and the next station that will seize the channel is also random. This type of protocols naturally leads to the formulation of polling models with random routing \cite{kleinrock, boxma}. A detailed comparison of various multi-access schemes, including random access schemes like ALOHA and CSMA, is given by Kleinrock \cite{kleinrock88}.

\paragraph{Optical networks.}
Polling models also find applications in the area of Ethernet Passive Optical networks (EPONs), where packets from different Optical Network Units (ONUs) share channel capacity in the upstream direction. An EPON is a point-to-multipoint network in the downstream direction and a multipoint-to-point network in the upstream direction. The Optical Line Terminal (OLT) resides in the local office, connecting the access network to the Internet. The OLT allocates the bandwidth to the Optical Network Units (ONUs) located at the customer premises, providing interfaces between the OLT and end-user network to send voice, video and data traffic. In an EPON the process of transmitting data downstream from the OLT to the ONUs is broadcast in variable-length packets according to the 802.3 protocol \cite{Kramer0}. However, in the upstream direction the ONUs share capacity, and various polling-based bandwidth allocation schemes can be implemented. Simple time-division multiple access (TDMA) schemes based on fixed time-slot assignment suffer from the lack of statistical multiplexing, making inefficient use of the available bandwidth, which raises the need for dynamic bandwidth allocation (DBA) schemes. A dynamic scheme that reduces the time-slot size when there are no data to transmit would allow excess bandwidth to be used by other ONUs. Kramer et al.\ \cite{Kramer1, Kramer2} propose an OLT-based interleaved polling scheme similar to hub-polling to support dynamic bandwidth allocation. To avoid monopolization of bandwidth usage of ONUs with high data volumes they propose an interleaved DBA scheme with a maximum transmission window size limit.

\paragraph{Bluetooth.}
Bluetooth is a wireless technology standard, used for exchanging data between mobile devices such as mobile phones, laptops, and headsets. These devices form small networks, referred to as Wireless Personal Area Networks (WPANs). The basic Bluetooth network topology is called a \emph{piconet}, consisting of one master device and up to seven slave devices. Miorandi et al. \cite{miorandi2004} observe that the structure of a piconet consisting of $N$ slaves, can be modeled adequately using a polling system consisting of $2N$ queues. One queue is used for each master-to-slave communication link, and one additional queue is required for each slave-to-master link. Approximations of the mean delays are found for the Pure Round-Robin (corresponding to $1$-limited service), gated, and exhaustive disciplines. Zussman et al. \cite{zussman2007} study the same model, deriving exact results regarding the packet delays.

\paragraph{I/O subsystems.}
Polling models occur in the context of I/O-subsystems of file servers or database management systems as well.
File servers are expected to handle large volumes of information-retrieval requests within a reasonable
amount of time. After the requested information (either static or dynamic) has been gathered, the
information is ready for transmission to the client over the communication network that connects the client
to the server. To this end, the information is placed in an I/O-subsystem, where it is ready to be drained
to the client over the network. This can be done for example via the Transport Control Protocol (TCP),
implementing a window-based control mechanism. The I/O-subsystem will typically consist of a number of
parallel TCP-buffers. The draining of the different buffers is managed by an I/O-controller that controls
the access of the buffer content to the network. To this end, the I/O controller ``visits'' the buffers in some order
to check whether a buffer has information to be drained, and if so, whether the corresponding TCP-window is open,
allowing for transmission of data segments residing in the buffer via the network to the client. We refer to
\cite{Reeser1, Reeser2} for performance models of Web servers including I/O-controlled buffers. This type of
I/O-control leads to the formulation of polling models where the service discipline represents the complicated
dynamics of TCP-window control. In \cite{czerniak}, Czerniak et al. analyze TCP systems using polling models to describe the different transmission methods.

\paragraph{Mobile networks.}
Polling models can also be found in the area of mobile networks, where different users compete for access to the shared
scarce radio resources. In such environments, the base station is typically in charge of assigning time slots to the different users in
some way. In this context, the server represents the right for transmission and the customers represent data packets to be transmitted.
Typical examples of polling mechanisms occur in the context of the Code Division Multiple Access (CDMA) based High Speed Packet
Access (HSPA), where the base station controller grants access to the medium in a per-timeslot basis. There are different scheduling
mechanisms for deciding which of the terminals gets access to the medium for the duration of a single time slot. A common implementation is simple Round-Robin (RR) scheduling, where medium-access is circulated among the terminals, independent of the quality of
the signal \cite{Berg1}. This immediately leads to polling models with limited service policies and cyclic server routing. A straightforward extension of RR scheduling is Weighted Round-Robin (WRR), which leads to the formulation of polling models with periodic server routing.
To enhance the efficiency of medium access for HSPA-networks, highly sophisticated channel-aware scheduling mechanisms have been proposed, granting time slots to terminals based on the measured instantaneous (possibly relative) Signal-to-Noise Ratios (SINRs) of each of the terminals
\cite{Bonald1, Borst22}. The inherent randomness in the channel conditions, and hence in the order in which the stations are granted
access to the medium naturally leads to the formulation of polling models with random or Markovian server routing
\cite{kleinrock, boxma}.

\paragraph{Ferry-based Wireless LANs.}
Polling models are applicable in the context of designing message ferry routes in Ferry-based Wireless LANs. In such FWLANs, a number of isolated nodes are scattered over some geographical area where communication between a node and the outer world, or communication between nodes, is made possible via a message ferry. The ferry follows a predetermined cyclic path, collecting messages from and delivering messages to nodes whenever it is \emph{in the vicinity} of the node (implying that the routes do not have to pass \emph{through} all the points in space).
The transmission and reception range is flexible: assuming a fixed transmission power, the range can be increased at the cost of decreasing throughput.
Taking into account the relation between the range and the throughput, one can design an optimal cyclic route of the ferry. We refer to Kavitha and Altman \cite{Altman_ITC} (and references therein) for results in FWLANs.

\paragraph{Mobile adhoc networks.}
Polling models occur naturally in the modeling of mobile ad-hoc networks (MANETs), consisting of both mobile and fixed wireless terminals. A typical feature of these kind of networks, is that wireless devices create their own wireless network in a distributional fashion. Mobile users can change location and thereby change communication links in the network. Examples of MANETs are animal-monitoring systems, collaborative conference computing, vehicular networks, peer-to-peer file-sharing systems and disaster-relief networks. Unlike most classical wireless networks, MANETs allow for multi-hop communication. Communication links appear and break down dynamically. Their lifetimes are random and independent of the amount of traffic offered to or transmitted over such links. To capture this phenomenon, De Haan \cite{Haan2} proposes the so-called pure time-limited service policy, where the server visits a queue for a random amount of time, independent of anything else in the system; note that this policy is not work conserving. The reader is referred to \cite{Haan1} (Chapter 1) for an overview of applications of polling models in MANETs.

\paragraph{Networks on chips.}
Another interesting application area of polling models is networks on chips (NOCs), which have been proposed as a solution for the inefficiency caused by traditional bus connections \cite{Benini2002, Dally2001}. In NOCs, intellectual property blocks (a general term for on-chip modules) are not connected to a single shared link, but to network interfaces that implement communication protocols. Data is transmitted using switches that consist of input and output ports.
If multiple ports have data for the same output port, only one input port can transmit its data, and the switch selects which one. Data that is not transmitted immediately is stored in buffers and will be transmitted later. We refer to \cite{Beekhuizen2010} for overview of the applicability and the state-of-the-art on analysis of polling models for networks on chips, and to \cite{vanhoudt2010} for an efficient numerical algorithm.

\subsection{Production systems}

A completely different application area of polling systems can be found in the  so-called \textit{stochastic economic lot scheduling problem} (SELSP). The SELSP deals with the make-to-stock production of multiple standardized products on a single machine with limited capacity under \textit{random demands}, \textit{possibly random setup times} and \textit{possibly random production times}. The SELSP is a common problem in practice, e.g., in glass and paper production, injection molding, metal stamping and semi-continuous chemical processes, but also in bulk production of consumer products such as detergents and beers. 

In many firms encountering the SELSP, the following class of \textit{fixed-sequence base-stock policies} is used for the control of the inventory of each product. We distinguish $N$ products, which are numbered $1,2,\ldots,N$. Subsequently, to each individual product a stock point is assigned which is controlled by a base-stock inventory policy. Under such a policy, for each product there exists a pre-defined desired number of items in stock, the base-stock level $b_i$, $i=1,2,\ldots,N$. When demand arrives at a stock point and the requested product is on stock, the demand is immediately fulfilled. Otherwise, demand is backlogged and fulfilled as soon as the product becomes available after production.
A production order, also called replenishment order, is placed immediately after demand for the corresponding product has arrived. These production orders queue up at the production facility, where each product has its own designated queue.

\begin{figure}[t]
\begin{center}
\includegraphics[width=1.0 \textwidth]{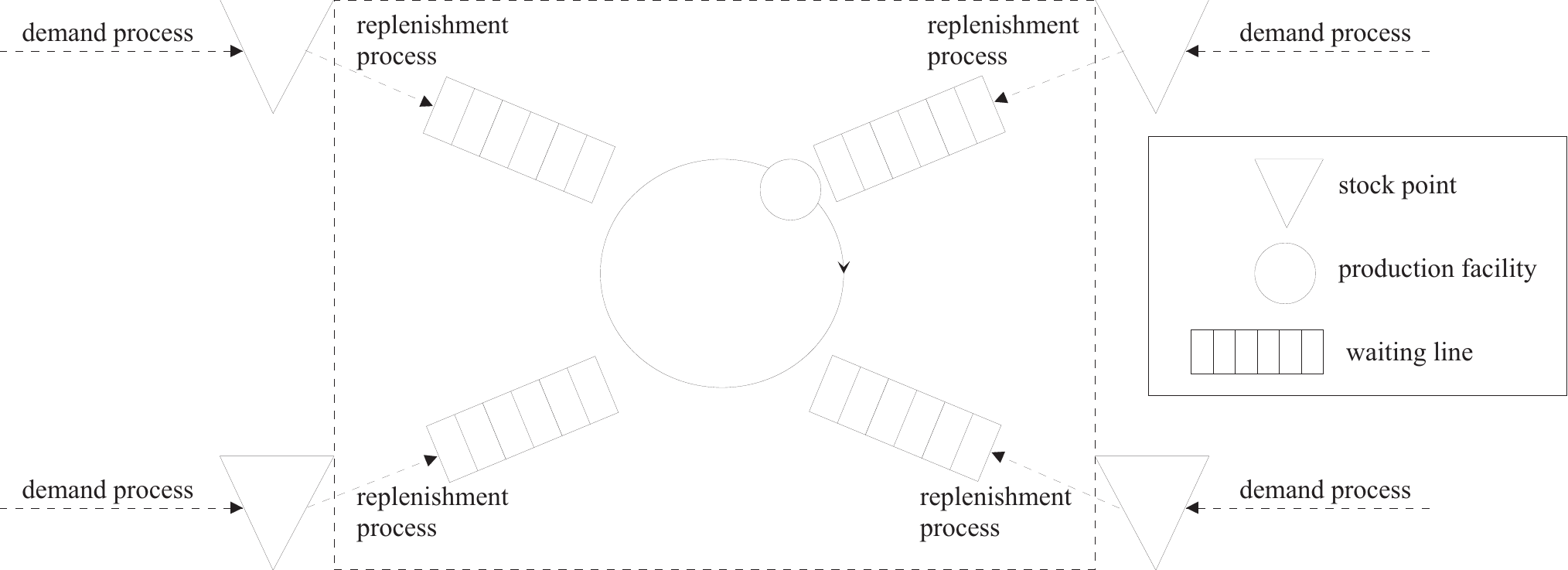} \\
\caption{The fixed-sequence base-stock system.}
\label{modelfig}
\end{center}
\end{figure}

On the strategy deployed by the production facility, the following two restrictions are imposed,
\begin{enumerate}
\item The products are produced according to a fixed production sequence;
\item When the machine starts production of a product, it will continue production until either the base-stock level has been reached or a second local criterion, i.e., only dependent on the stock level of the product currently setup, has been fulfilled.
\end{enumerate}
The combination of stock points and production facility is visualized in Figure \ref{modelfig}.

For given values of the base-stock levels $b_i$ the steady-state net stock level $N_i$ for product $i$ is given by, for $i=1,2,\ldots,N$,
\begin{equation}
N_i=b_i-L_i,
\end{equation}
where $L_i$ denotes the steady-state \textit{shortfall} (the number of outstanding production orders at the production facility) of product $i$. Notice that the net stock level of a product becomes negative, when the shortfall of this product is larger than its base-stock level. One can verify that the shortfall of a product is independent of the base-stock levels implying that the performance of the production facility can be analyzed independently of these base-stock levels. Moreover, the shortfall distribution of a product at the production facility is identical to the queue length distribution in a polling system. The interarrival, service and setup time processes in such a polling system are identical to the demand, processing and setup time processes in the SELSP, respectively.
\textit{For given base-stock levels, the evaluation of a fixed-sequence base-stock policy is, therefore, tantamount to evaluation of the corresponding polling system}.

\paragraph{Literature on the SELSP.} Seminal papers analyzing the above fixed-sequence base-stock policy via polling systems are by Federgruen and Katalan \cite{Federgruen2,Federgruen3,Federgruen5}. Besides the basic assumptions for the SELSP, 
products are produced by an exhaustive or a gated base-stock policy. The production manager is allowed to insert a fixed idle time prior to the setup for a product in order to reduce the setup frequencies and, hence, the average setup costs. 

Grasman et al. \cite{Grasman1} extend the exhaustive base-stock model of Federgruen and Katalan by adding random yields for the cases of backlogging, lost sales and expediting. In case of backlogging they derive similar newsboy equations in order to obtain optimal base-stock levels. In case of lost sales or expediting they have to resort to a heuristic for computing the (approximate) optimal base-stock levels. Krieg and Kuhn \cite{Krieg1,Krieg2} introduce continuous-time models for single-stage multi-product \textit{Kanban systems}, which are completely identical to the SELSP with lost sales. 
In Krieg and Kuhn \cite{Krieg1}, a system is analyzed with \textit{state-independent} setups, whereas in Krieg and Kuhn \cite{Krieg2} \textit{state-dependent} setups are modeled, i.e., no setup for a product is incurred when there is no shortfall. Production quantities are in both models determined by an exhaustive base-stock policy. 

Winands et al. \cite{winands2} analyze a two-product system, in which a high-priority product is produced exhaustively and a low-priority product according to the $k$-limited service strategy. 
Chapter 4 of \cite{winands} is devoted to a numerical simulation study which shows that the $k$-limited lot-sizing policy outperforms the standard exhaustive policy for a wide variety of environments. In particular, the $k$-limited policy proves its worth in asymmetric production systems.

In most applications, the service discipline determines how many items are produced, and, hence, what the (random) duration of a cycle in the polling system is. In some cases, however, the cycle length is fixed implying that the corresponding analysis often reduces to a one-dimensional problem. See \cite{giezenaar} for a case study of a chemical plant, for which a fixed production sequence strategy in combination with a fixed cycle length has been developed under the assumption of deterministic production and setup times. Erkip et al. \cite{erkip}  introduce a discrete-time model under the assumption of backlogging, in which the production and setup times are deterministic. Other work in this direction is by Bruin \cite{bruin}, who presents a generating function approach for the fixed cycle strategy under general traffic settings, and by Dellaert \cite{dellaert}, who develops a heavy-traffic approximation for the optimal base-stock levels.

\paragraph{Make-to-order.}
Summarizing, we can say that the SELSP is an extension of a standard polling system by an additional inventory dimension.
The make-to-order scheduling counterpart (in which no stock is kept for a product and thus the base-stock levels equal 0) is, however, precisely equivalent to a standard polling model. The polling system again consists of a server, mostly referred to as machine, and
multiple customer classes, which will be referred to as products.

In the literature, one can notice that (almost) all modeling variants as
discussed in Section \ref{pollingsection} have been studied in the context of such a make-to-order
production situation. This implies that an enormous amount of papers on polling
systems have appeared motivated by make-to-order manufacturing applications.
Since the primary goal of the present paper is to give an overview of the
various application areas of polling systems - and not to give an exhaustive
list of papers within each area - we only give a small subset of recent papers
on make-to-order polling systems (i.e., \cite{maketoorder7,maketoorder2,maketoorder3,maketoorder1,maketoorder8,maketoorder4,maketoorder6,maketoorder5}).

Finally, in Markowitz et al. \cite{Markowitz1,Markowitz2}
heavy traffic polling results are applied to all kinds of related stochastic
multi-product single-machine scheduling problems.

\paragraph{The impact of setup times.}
Besides papers on the performance analysis of polling systems in the context of production environments, much research has been published on the impact of setup times (see \cite{Duenyas1,Gerchak1,McIntyre1,Sarkar1,Zangwill1,Zangwill2}). In these papers, it is shown that reduction of setup times can, counterintuitively, increase the mean queue lengths in polling systems and cyclic production systems for a variety of settings. Furthermore, the conditions are studied under which mean queue lengths increase when setup times are reduced \cite{Samadder2007,Samadder2008,Federgruen1,Federgruen2,winands}.

\subsection{Traffic and transportation systems}

The third important area in which polling systems are frequently applied in practice are traffic and transportation systems. In road traffic a polling system is the natural way to model a situation where queues arise due to the fact that multiple flows of traffic have to share one single lane. The plainest form is a two-way road which is partly blocked because of an accident or road maintenance, but any common traffic intersection also qualifies. Besides road traffic, polling models are employed in transportation systems
consisting of driverless electric vehicles that follow a predetermined track. The actual transportation is performed by so-called Automated Guided Vehicles (AGVs). AGV systems are used in warehouses and container terminals, but also in other application areas like, e.g., health care.

\subsubsection*{Traffic signals}

We discuss applications of polling systems in road traffic first. Queues are formed by the lines of cars waiting before a red traffic signal. The time that is required for one vehicle to pass the traffic signal can be viewed as service time, while the times that all queues face a red light simultaneously, i.e. the clearance times of the intersection, can be considered as switch-over times. There are several aspects that make polling models for signalized intersections different from the models in other application areas.
Firstly, the assumption of independent, identically distributed service times is in practice not valid. When a green period commences, a certain time elapses while vehicles are accelerating to normal speed. After this time, which typically lasts a few seconds, the queue discharges at a more or less constant rate, which is called the saturation flow. In practice the saturation flow may vary within a cycle and between cycles, but Webster \cite{webster58} finds that the assumption of a constant saturation flow agrees well with values observed over a fairly large number of cycles. A typical feature of the traffic light queue is that cars approaching an empty intersection while facing a green light do not slow down and require almost no service time at all, especially if they do not have to take a turn.
The second difference between traffic intersections and many other polling applications is the arrival process. 
The assumption of Poisson arrivals might be realistic for an isolated traffic intersection, but many intersections are part of an arterial system, which means that an output process of the first intersection forms an input process for the next intersection. Only few papers deal with these so-called \emph{platooned} arrivals.
The third difference that we discuss is that most traffic intersections are divided into groups of flows that face a green light simultaneously. In the polling system this can be modeled as a (possibly varying) number of multiple coupled servers serving queues simultaneously. The moment at which a switch-over is incurred depends on the service discipline that is used for the intersection. The best known discipline is the \emph{Fixed Cycle Traffic Light} discipline, which uses deterministic green, red and amber times. The obvious disadvantage of fixed settings is that the system does not respond to the current situation, which might be less efficient because traffic lights remain green even when no cars are waiting in the corresponding queue. An alternative to fixed settings that is very commonly used is an \emph{adaptive control mechanism} for traffic signals. This is a mechanism that detects the presence of vehicles and makes decisions on whether or not to switch a signal based upon this information. A typical policy is to wait until all traffic flows that face a green light are empty before switching to red. 
For some intersections, a time limited policy might be more suitable. See \cite{newell98} for a discussion on this topic. We discuss literature on polling models for traffic intersections in more detail by splitting it into two groups: papers concerning the fixed cycle traffic light queue, and papers on vehicle-actuated signals.

\paragraph{Literature on the Fixed Cycle Traffic Light queue.}
A traffic intersection with fixed green, amber and red times can be viewed as a polling model with deterministic visit times and switch-over times. One may argue whether this system should be considered a real polling system, since the analysis of delays is a one-dimensional problem and one can focus on one queue in particular. One of the most influential papers on the fixed cycle traffic light queue is \cite{webster58}, in which Webster describes a procedure to find optimal settings (i.e., green and red times for all traffic flows) for traffic lights. He derives an approximate expression for the optimal cycle length and an expression for the average delay per vehicle. His expressions are partly based on theoretical grounds, and partly on simulation results. Although more sophisticated methods have been developed ever since, his results are still quite popular in practice. Few years later, Miller \cite{miller63}, under the assumption of Poisson arrivals, and Newell \cite{newell65}, for general arrivals, develop approximations for the distributions of the delay and the queue length. Exact expressions for these distributions have been obtained by Heidemann \cite{heidemann94}, but again under the assumption of Poisson arrivals.
In a recent paper Van Leeuwaarden \cite{leeuwaarden06} derives the queue length and the waiting time distribution for general arrival processes. Van den Broek et al. \cite{vandenbroek06} derive approximations for the mean overflow, i.e., the mean queue length at the end of a green period, that are easier to compute and provide more intuitive insight.

In practice, intersections are often part of an arterial system and interarrival times are correlated. Alfa and Neuts \cite{alfaneuts95} regard a model with platooned arrivals, which allows for distinguishing between interarrival times between cars within a platoon, and interplatoon interarrival times. The number of cars within one platoon can have any discrete phase probability distribution. The platooned arrival process is applied to road traffic, and to the fixed cycle traffic light queue. Through a numerical example they conclude that ignoring correlation in the arrival process leads to an underestimation of the mean queue length at high traffic intensities.


\paragraph{Literature on vehicle-actuated traffic signal control.}

Nowadays fixed cycles are less and less frequently implemented at traffic intersections, but it might still be a realistic assumption when the intersection is congested during rush-hour. Most traffic signal systems are vehicle-actuated, which means that the system contains detectors that gather information about the number of cars present at each flow and regulate the green and red times based on this information. In particular, when a queue becomes empty, the corresponding traffic light turns red and the traffic light of the next group of flows becomes green. This situation corresponds to a polling model with exhaustive service. However, in general, vehicle-actuated systems are designed to have minimum and maximum green times for each flow. This makes them very suitable to model as a polling system with time-limited service, but it also makes them difficult to analyse. 
Frigui and Alfa \cite{friguialfa98,friguialfa99} propose an iterative algorithm to approximate waiting times in a discrete-time polling system with time-limited service. Their papers focus on applications in communication systems, but in \cite{friguiPhD} the model is also applied to traffic intersections.

The earliest literature on vehicle-actuated systems dates from the early 1960s when Darroch et al. \cite{darroch64} analyze a system that consists of two intersecting traffic streams that are served exhaustively. A model with two lanes that does not assume Poisson input has been studied by Lehoczky \cite{lehoczky72}, who uses an alternating priority queueing model. Newell \cite{newell1} analyzes an intersection with two one-way streets  using fluid and diffusion queueing approximations.
Daganzo \cite{Daganzo} studies a polling system with more general arrival and service processes, with applications in traffic and transportation systems (as long as only one flow of traffic is served at a time).
In \cite{newell2}, Newell and Osuna study a four-lane intersection where two opposite flows face a green light simultaneously.
A variation of the two-lane intersection is introduced by Greenberg et al. \cite{greenberg88}, who analyse mean delays on a single rail line that has to be shared by trains arriving from opposite directions. This model is extended by Yamashita et al. \cite{yamashita06}, who study alternating traffic crossing a narrow one-lane bridge on a two-lane road.
In many of the discussed papers traffic is modeled as fluid passing through the road. This approximation is fairly accurate when the traffic intensity is relatively high. Vlasiou and Yechiali \cite{vlasiouyechiali07} use a different approach, modeling a traffic intersection as a polling system with an infinite number of servers visiting each queue simultaneously.

All in all it is rather surprising how few mathematical models have been developed that deal with intersections where multiple flows face a green light simultaneously. In \cite{boon} a two-queue polling model is analyzed with two priority levels in one queue. This model can be applied to an intersection where two \emph{conflicting} traffic flows share the same green slot.
In \cite{haijemavanderwal07} a Markov Decision Problem (MDP) decomposition approach is used to find optimal traffic signal settings at intersections with larger numbers of combined flows.

\subsubsection*{Automated guided vehicles}

Next to traffic intersections, a large transportation related application area of polling systems are Automated Guided Vehicles. The first AGV system was simply a tow truck that followed a wire in the floor instead of a rail, nowadays AGVs are mostly laser navigated. AGVs are mostly used for the transport of materials in manufacturing systems, but also in other areas like public transportation systems at airports. Not all AGV systems can be modeled as a polling system, but we discuss some exceptions.  In a conventional AGV system, each vehicle can pick up a load from any station and deliver it to any other station. In order to avoid collisions between vehicles, most systems use the \emph{zone blocking} concept where the entire system is divided into zones. The control system allows only one vehicle in each zone at a time. Whenever an AGV system consists of a single loop, it can be modeled as a polling system. The vehicle corresponds to the server in the polling system, and the stations form the queues of the system. In some cases each station is modeled as two queues, one for the picking-up, and one for the drop-off. The inter-station travel times are modelled as switch-over times in the polling model. AGV systems with multiple vehicles can be modeled as a polling system with multiple, independently moving servers, but it is more common to approximate this situation by modeling it as a system with one vehicle that moves at a higher speed (see, e.g., \cite{dukhovnyy}). It is typical for AGV systems that the vehicle in general does not visit the stations in a deterministic order. Although most papers concerning polling systems focus on a fixed visiting order of the queues, several papers discuss an alternative server routing. Srinivasan \cite{srinivasan91} uses a polling system with Markovian server routing to model AGV systems. In \cite{boxmaoptimization} and \cite{bertsimasxu} optimization of server routing in polling systems is discussed. In AGV systems the visit order is generally determined by a decentralized cargo dispatching rule. The most commonly used cargo dispatching rule in AGV systems is First-Encountered-First-Served (FEFS) rule. The FCFS rule is less efficient, because it leads to unnecessary empty travel. Note that FCFS in the context of AGV systems generally does not refer to the order in which items are picked up within a queue, but to the route that the vehicle takes. In terms of queues, it travels towards the queue with the oldest waiting customer. The FEFS rule is presented for single-loop AGV systems in \cite{bartholdi} and states that an empty vehicle continues to travel along the loop until it finds some load to pick up at some station. If it has available space, it picks up the load and drops it off at its destination. The number of packages that are picked up depends on the service discipline and may vary per station.

\begin{figure}[t]
\includegraphics[width=0.47\linewidth]{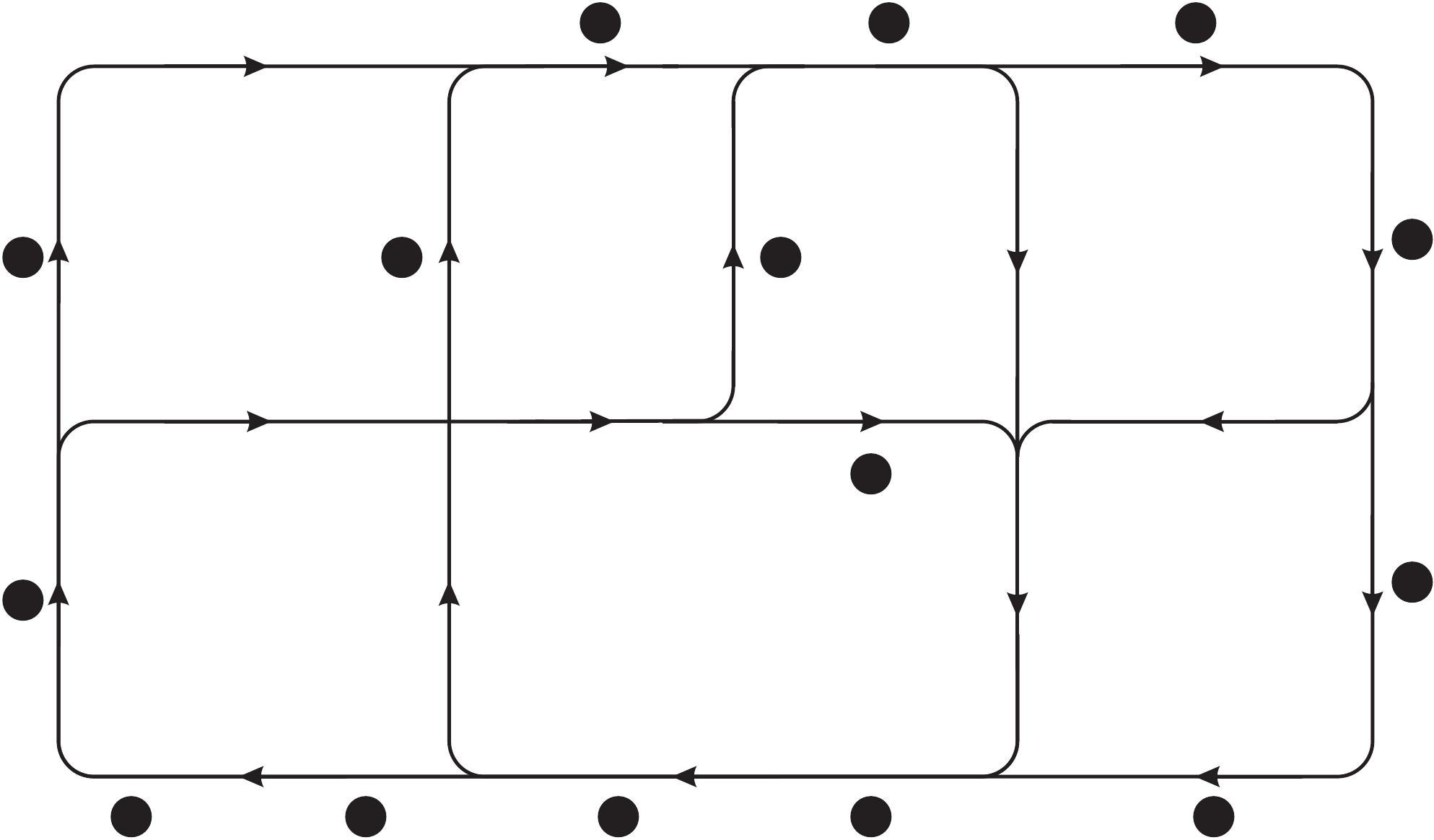}
\hfill
\includegraphics[width=0.47\linewidth]{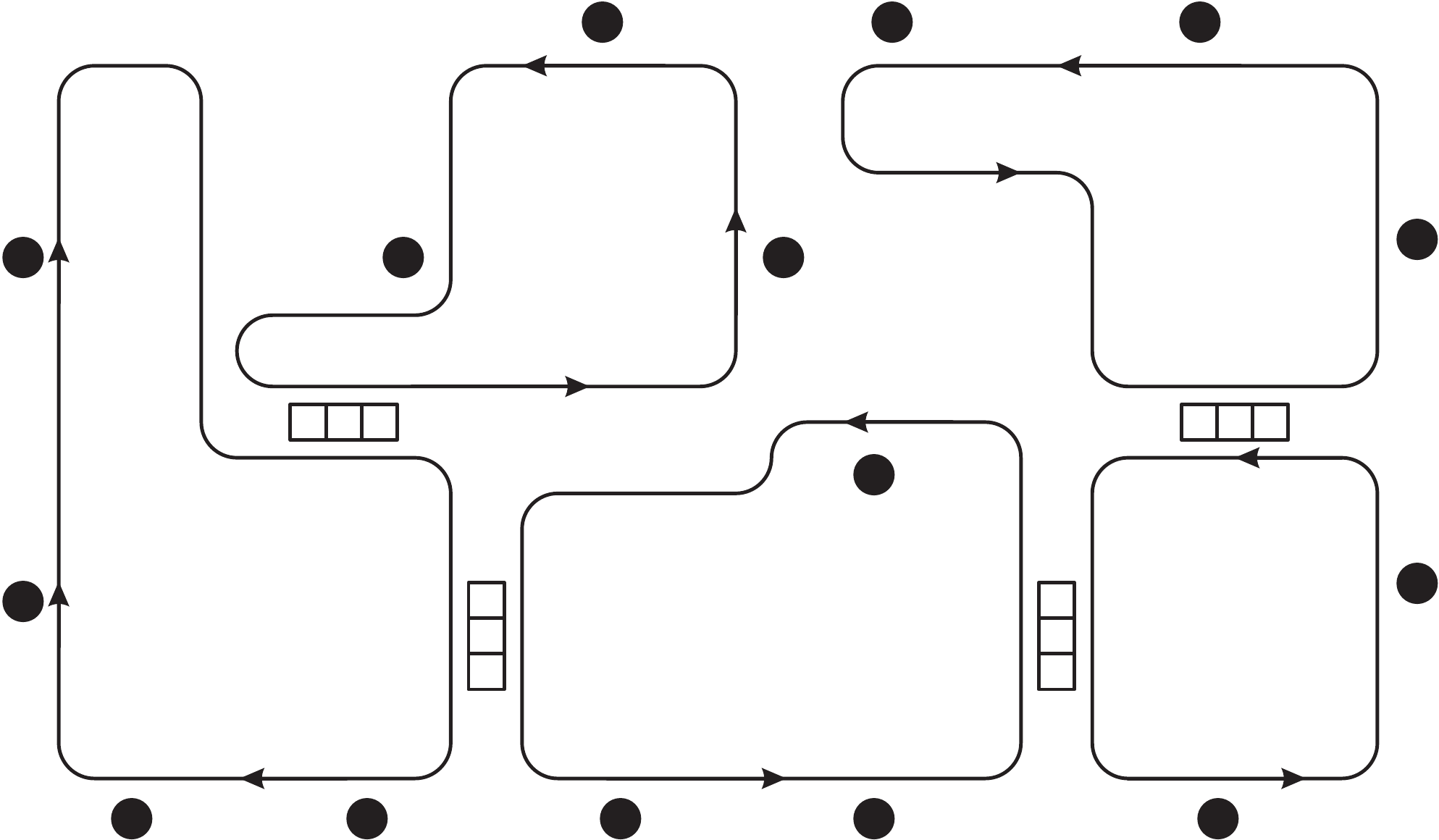}
\caption{A conventional AGV system (left) and a tandem AGV system (right) as proposed by Bozer and Srinivasan \cite{bozersrinivasan91}. The dots represent stations with pick-and-drop points. In the tandem AGV system, each loop is modeled as a single server polling system, and additional pick-and-drop points are placed as interface between adjacent loops. This figure appeared originally in \cite{bozersrinivasan91}.}
\label{agvpicture}
\end{figure}

Bozer and Srinivasan \cite{bozersrinivasan91} discuss an AGV system with tandem configurations. The original AGV system is divided into non-overlapping, single vehicle closed loops with load transfer stations in between (see Figure \ref{agvpicture}). Each loop can be modeled as a single server polling system. An advantage of a tandem AGV system is that each vehicle is dispatched over a smaller number of stations. Another advantage is that traffic management problems within each loop are completely eliminated. In \cite{bozersrinivasan91} the mean throughput capacity for the loop is estimated under the FEFS rule. Srinivasan et al. \cite{srinivasanbozer94} study this rule as well, and also introduce a modified FCFS dispatching rule, under which a vehicle deposits a load at a station and then checks that station for a new load-transfer request. If there is one, the vehicle serves this request; otherwise, it travels empty to other stations based on the FCFS dispatching rule. Xu et al. \cite{xushicheung07} focus on the loop configuration rather than on the cargo dispatching rule. Ganesharajah et al. \cite{ganesharajah} consider both simultaneously. Van der Heijden et al. \cite{vdheijden01} study an underground transportation system that contains a single traffic lane which has to be shared by AGVs from opposite directions. The system has similarities with a road traffic situation, but the extremely long clearance times require a more dedicated approach. The research is continued in \cite{ebben04} where more intelligent, adaptive control rules and dynamic programming algorithms are used to minimize vehicle waiting times.

\subsection{Miscellaneous}

The present subsection provides an overview of miscellaneous applications of polling systems observed
in the open literature.

\paragraph{Health care.}

In \cite{cicin} a medical emergency room is modeled via a multi-server polling model. At the emergency
room different types of patients arrive and are placed in queues
depending on the type of surgical procedure required. Each emergency
surgical procedure has a separate queue of infinite buffer capacity.  The number of emergency room theatres is limited, and each emergency room theatre is modeled by a server in the polling model. Since setting up for surgery procedures takes longer than the actual
procedures themselves, the service times within the polling model are much smaller than the switch-over times.
Finally, an urgency parameter in terms of average waiting time is assigned to each patient which dictates the existence of (local) priority levels within each queue. From the application point of view, Cicin-Sain offers, albeit in a highly idealized mathematical setting, a preliminary exploration of how to design emergency rooms such that the most urgent patients have the shortest waiting times  \cite{cicin}. Another application in health care is depicted in \cite{vlasiouadanboxma}, where a polling model, consisting of two queues with infinite supply, is used to model scheduling strategies by surgeons performing eye surgeries.

\paragraph{Mail delivery.}

Another polling application is an internal mail delivery system, which was
modeled via a continuous polling system (a system with an infinite number of queues) by Nahmias and Rothkopf \cite{nahmias}. In their model,
a clerk traverses at a constant rate a cyclic route along which mail is generated according to a Poisson process.
The clerk picks up the mail that has been generated since the last traversal and delivers the mail (that was previously picked up and sorted) to locations distributed uniformly along the route. At the end of the route, the clerk sorts the mail just picked up. Then, the clerk again traverses the route, delivers the mail that has just been sorted and picks up the new mail and so on. In \cite{nahmias} also the extension is studied in which there are several independent routes for a couple of clerks and a single room where sorting takes place.

Sarkar and Zangwill \cite{Sarkar2} discuss a finite-queue variant of the above problem. That is, they study a polling system where the workload at a queue either comes from outside the system or from another queue within the system. In terms of the above mail delivery application the sorting of the clerk would be done at one designated queue $N$, whereas at the other queues the mail is generated by exogenous Poisson processes. The sorting work done by the clerk at queue $N$ obviously depends upon the amount of mail generated at the other queues. Besides the mail delivery application, Sarkar and Zangwill also describe applications of this model related to rework in manufacturing systems, computer file transfers and buses circulating at airports.

\paragraph{Snow blowing.}

Eliazar \cite{eliazar1,eliazar2} studies a so-called snowblower problem. That is, on a closed-loop
racetrack snow is falling randomly and a snowblower machine is continuously
circling this racetrack and clearing off the snow.  This snowy racetrack can be modeled
by a continuous polling system  by drawing
the following analogies: server $\leftrightarrow$ snowblower; job arrivals $\leftrightarrow$ snowfall; workload
$\leftrightarrow$ snowload.

\paragraph{Shipyard loading.} Another application of a  polling system is a shipyard loading problem in which containers arrive by truck to the port. Trucks drive the containers up to their destination and the receiving yard crane
transfers them to their assigned position. The crane represents the server of a polling system, whereas the trucks arriving for a specific destination are customers of a specific type. As stated by Daganzo \cite{Daganzo}, storage room (or: queue lengths) is the most critical performance measure for transportation problems dealing with freight. 
For transportation problems dealing with passengers, the waiting time is the most commonly used performance measure. In the aforementioned shipyard loading problem, one is especially interested in
 the characteristics of the \textit{total} queue length in the polling system. That is, early arriving trucks
wait in a common area until they are due for service. To dimension up this common area, information on the total queue length is needed.


\paragraph{Dynamic Picking Systems.} Gong and De Koster \cite{gongdekoster08} use a polling model to describe a dynamic order picking system (DPS). In a DPS, a worker picks orders that arrive in real
time during the picking operations and the picking information can dynamically change in a picking cycle.
It is very important to achieve short delivery times, especially for online retailers. One of the challenging questions that online retailers now face is how to organize the logistic fulfillment processes during and after order receipt. In traditional stores, purchased products can be taken home immediately. However, in the case of online retailers, the customer must wait for the shipment to arrive.
In \cite{gongdekoster08} polling models are used to describe and identify a DPS for online retailers. Polling-based picking systems can lead to shorter throughput times than traditional batch picking systems, particularly for high order arrival rates.

\paragraph{Elevators.}

A (multi-server) polling model has been studied by Gamse and Newell \cite{gamse1,gamse2} for an application
of elevators in a building. In the model assumptions are made pertaining to the relative movements of servers within the system, with the assumptions being formed based on the nature of the specific application of elevators. In this context it is also interesting to spend some words on the so-called \textit{elevator} server routing scheme. In such a system, the server first serves queues in the ``up'' direction, i.e., in the order $1,2,\ldots,N-1,N$, and subsequently serves these queues in the opposite (``down'') direction, i.e., visiting them in the order $N,N-1,\ldots,2,1$ (see, e.g., Altman et al. \cite{altman5}). It turns out that the globally gated service discipline in combination with elevator-type server movement results in equal mean waiting times for customers in \emph{all} queues, establishing a perfectly fair polling system.

Another way of modeling an elevator system using a polling model can be found in \cite{gandhi96}. In this paper an MDP approach is used to dynamically schedule elevators in a building. However, the model under consideration is only suitable for a limited amount of practical purposes, because all customers are assumed to have one common destination, which is the ground floor.

Besides the different routing scheme, another difficulty in modeling elevator systems is the fact that each idle elevator should return to the floor that is its home base. In polling literature some work has been done on this topic, which is called a stopping or dormant server, cf., \cite{borst,eisenberg94}.

\paragraph{Maintenance.}

We would like to end this survey with actually the first application of a polling system which appeared in the open
literature, i.e., in the area of maintenance. That is, Mack et al. \cite{mack2,mack1} use a polling
model to describe a patrolling repairman who inspects a number of machines to check whether a breakdown has
occurred and if so, eliminates such breakdowns. Evidently, in a polling model the repairman is represented by
the server, the breakdowns are represented by the customers and the times needed by the repairman to travel
from one machine to the next are represented by the switch-over times. In \cite{koenigsberg} a similar model is
studied in which an operator at a fixed position serves a number of storage locations on a rotating carousel conveyor.
Models with several independent rotating carousels have also been considered, cf. \cite{bunday,kim}. Weststrate \cite{weststrate}
studies a polling model which captures the behavior of a single repairman who is not only concerned with corrective maintenance,
i.e., maintenance after a breakdown has occurred, but who can also perform preventive maintenance.

\section{Outlook}

Let us conclude the present survey with some personal views on possible future research
directions related to applications and practical considerations. The present section is structured identically to the preceding section.

\subsection{Communication systems}

The technological advances in information and communication technology lead to the formulation of polling models with
features that are not included in the classical polling models framework for which results are known today, and as such
will lead to a variety of scientific challenges in the years to come. In this subsection we address a number of such open
research problems.

The applicability of the results on polling models that are known today is strongly limited by the usual {\it Poisson assumption}, often needed to obtain exact results. Despite the fact that the packet-level dynamics on communication networks are well-known to be highly complicated and often far from Poisson, the vast majority of research results explicitly rely on the Poisson assumptions, or some straightforward relaxation of that. Extension of Poisson-based results to non-Poisson and/or dependent arrival processes raises tremendous scientific challenges, requiring new queueing-theoretical analysis methods. The development of such methods would strongly enhance the applicability of queueing-theoretical results on polling models. 
In this context, encouraging results (but still only for renewal arrivals) have been obtained in \cite{Bertsimas1,mei1,vdmeiwinands08,olsenvdmei05}.

Another aspect that limits the applicability of results on polling models is {\it inflexibility with respect to the service
policies}. In fact, the modeling of many computer-communication networks often leads to polling models where the service policies, typically
representing the right for transmission/receipt of data to a shared medium, are of limited-type, and as such do not possess the
convenient branching structure that allows for exact analysis. This raised the great need for the development of
techniques to handle limited-type service policies, although there is not much hope for detailed exact analysis of such models.
This stresses the importance of the development of asymptotic results (e.g., asymptotics on heavy traffic, large numbers of queues,
heavy tails), a direction of research that has been highly successful in analyzing many queueing models in the past couple of years.
Development of asymptotic results on polling models with limited-type service policies will continue to form a challenging area
of research in the years to come.

The assumption of having a \emph{single} server serving multiple queues is also an important limitation for many real applications. As mentioned before, only few results exist on polling models consisting of multiple servers.

In the area of mobile and wireless networks, a promising means to increase bandwidth, and to enhance robustness of the user-perceived
service quality, is to make use of multiple antennas, enabling users to simultaneously utilize multiple networks (e.g., MIMO). In such
environments, packet streams can be split over multiple parallel networks, each of which can be modeled as a polling model,
where the packet routing may depend on the actual state of these networks. The packet-level analysis of such multi-network scenarios
using this type of concurrent-access techniques will lead to the development of coupled polling models, where the arrival processes at
the different polling models are correlated and state dependent. For this type of models, only few results are known today. Possibly, results from
\cite{levy5} (correlated arrivals) and \cite{pollinglevy09} (correlated and state dependent arrivals) could be used to model these kind of networks.

\subsection{Production systems}

Although polling models appear naturally as a modeling tool in a wide variety of production applications,
one can certainly observe common characteristics among these applications.
The first similarity is the fact that a \textit{high utilization of capacity} is typically prevalent
due to either the magnitude of the setup times or the customer demands. 

Secondly, in \cite{inman,winands3} it is shown that assuming exponential interarrival times is not appropriate, whereas
the independence assumptions do appear to be valid. Although these studies represent only a minuscule sample of production systems,
they do show that the analysis of queueing systems with \textit{general (renewal) arrival processes}
is of practical interest for production applications.

Thirdly, in many other production systems the single most important performance measure is not an aggregate measure like
the mean queue length, rather the probability that the queue length exceeds a pre-defined threshold.

In view of these three observations, the importance of explicit expressions of the complete queue length distribution in the asymptotic regimes of high utilization of capacity (due to either the presence of setup times
or the customer demands) for general arrival processes is evident. 

\subsection{Traffic and transportation systems}

Although traffic intersections with vehicle actuated signals are typical examples of polling systems, there is still a lot of research to be done in this area before a realistic model is obtained. In traffic and transportation systems the assumption of independent interarrival times is reasonable for isolated intersections, but not for intersections that are part of an arterial system. Just like in the area of communication and production systems, it is necessary to develop techniques to analyse polling systems with other arrival assumptions, in particular platooned arrivals. For practical purposes there is a great need of good models for intersections that are part of an arterial system, i.e., the output of a certain queue forms the input process of a queue in a neighboring polling system. The assumption of independent interarrival times is strongly violated in such systems. One could consider the study of Reiman and Wein \cite{reimanwein99} as an important first step in this direction. They study waiting times of customers in a sequence of two-queue, exhaustive polling systems under heavy traffic conditions. Obviously, the analysis of networks of polling systems would not only be of interest for traffic and transportation systems, but also for other application areas.

The assumption of independent service times is not realistic either in most cases. It is also of great importance that models are developed for situations where multiple (conflicting and/or non-conflicting) traffic flows face a green light simultaneously. Some initial results on this topic are obtained in \cite{boontrafficlights2010,boonthesis}, but (also here) many challenges remain for future research. Finally, and also this is not just important in the area of transportation, it is important to obtain good waiting-time approximations for time-limited systems.

Another area in transportation systems where few results have been obtained, is a system with one or multiple elevators. It would be of much practical relevance to have a realistic (polling) model that can be used to find optimal settings for the visiting policies of elevators. Difficulties that have to be addressed in elevator polling are the different routing schemes, but also the fact that each idle elevator returns to the floor that is its home base.

\section*{Acknowledgments}

The authors are very grateful to Ivo Adan and Onno Boxma for providing valuable comments on earlier drafts of the present paper.

\small
\bibliographystyle{marko}

\end{document}